\documentclass{article}

\usepackage{amsmath}   
\usepackage{amssymb}

\newtheorem{Th}{Theorem}  
\newtheorem{Prop}{Proposition}   

\newtheorem{Lemma}{Lemma}   
\newtheorem{Coro}{Corollary}

\newcommand{\finishproof}{\hfill $\Box$ \vspace{3mm}}
\newcommand{\R}{\mathbb{R}}
\newcommand{\Z}{\mathbb{Z}}
\newcommand{\C}{\mathbb{C}}
\newcommand{\T}{\mathbb{T}}
\newcommand{\h}{{\mathfrak h}}

\newcommand{\al}{\alpha}
\newcommand{\ga}{\gamma}
\newcommand{\la}{\lambda}
\newcommand{\Om}{\Omega}

\begin{document}

\date{}

\title{On the symplectic phase space of KdV}   
 
\author{T. Kappeler\footnote{Supported in part by the Swiss National Science Foundation, and the programme
SPECT, and the European Community through the FP6 Marie Curie RTN ENIGMA
(MRTN-CT-2004-5652).}, F. Serier, and P. Topalov} 
  
\maketitle

\begin{abstract}  
We prove that the Birkhoff map $\Om$ for KdV constructed on $H^{-1}_0(\T)$ can be interpolated between
$H^{-1}_0(\T)$ and $L^2_0(\T)$. In particular, the symplectic phase space $H^{1/2}_0(\T)$ can be described
in terms of Birkhoff coordinates.
\end{abstract}

\section{Introduction}
In \cite{KT3} it is shown that the Birkhoff map for the Korteweg - de Vries equation (KdV),
on the circle $\T:=\R/\Z$, introduced and studied in detail in \cite{KP,KM} can be analytically extended to
an analytic diffeomorphism
\[
\Om : H^{-1}_0(\T)\to\h^{-1/2}
\]
from the Sobolev space of distributions $H^{-1}_0(\T)$ (dual of $H^1_0(\T)$) to the Hilbert space
of sequences $\h^{-1/2}$ where for any $\al\in\R$,
\[
\h^\al:=\{z=(x_k,y_k)_{k\ge 1}\;|\;\|z\|_\al<\infty\},
\]
with
\[
\|z\|_\al:=\Big(\sum_{k\ge 1}k^{2\al}(x_k^2+y_k^2)\Big)^{1/2}\,.
\]
In this paper we show that $\Om$ can be interpolated between $H^{-1}_0(\T)$ and $L^2_0(\T)$.
\begin{Th}\hspace{-2mm}{\bf .}\label{Th:main}
For any $-1\le\al\le 0$,
\[
\Om|_{H^\al_0(\T)} : H^\al_0(\T)\to\h^{\al+1/2}
\]
is a real analytic diffeomorphism.
\end{Th}
As an application of Theorem \ref{Th:main} we characterize the regularity of a potential
$q\in H^{-1}(\T)$ in terms of the decay of the gap lengths $(\ga_k)_{k\ge 1}$ of the periodic spectrum 
of Hill's operator $-\frac{d^2}{dx^2}+q$ on the interval $[0,2]$. More precisely,
recall that the periodic spectrum of $-\frac{d^2}{x^2}+q$ on the interval $[0,2]$ is discrete. 
When listed in increasing order (with multiplicities) the eigenvalues $(\la_k)_{k\ge 0}$ satisfy 
\[
\la_0<\la_1\le\la_2<\la_3\le\la_4<...
\]
The gap lengths $\ga_k=\ga_k(q)$ are then defined by
\[
\ga_k:=\la_{2k}-\la_{2k-1}\,\,\,(k\ge 1)\,.
\]
\begin{Th}\hspace{-2mm}{\bf .}\label{Th:characterization}
For any $q\in H^{-1}(\T)$ and any $-1\le\al\le 0$, the potential
$q$ is in $H^\al(\T)$ if and only if $(\ga_k(q))_{k\ge 1}\in\h^\al$.
\end{Th}
In a subsequent paper we will use Theorem \ref{Th:main} to study the solutions of the KdV equation
(see \cite{CKSTT,CKSTT1}, \cite{KT5}, \cite{TaoBook}) in the symplectic phase space $H^{-1/2}_0(\T)$ introduced by
Kuksin \cite{Ku}.

\vspace{0.3cm}

\noindent{\em Method of proof:} Theorem \ref{Th:characterization} can be shown to be a consequence of
Theorem \ref{Th:main} and formulas relating the $n$'th action variable $I_n$ with the $n$'th gap length $\ga_n$
and their asymptotics as $n\to\infty$. In view of results established in \cite{KT3} the proof of Theorem \ref{Th:main}
consists in showing that for any $-1<\al<0$ the restriction of $\Om$ to $H^\al_0(\T)$,
$\Om|_{H^\al_0(\T)} : H^\al_0(\T)\to\h^{\al+1/2}$, is {\em onto}. Our method of proof combines a study of the Birkhoff map
at the origin together with a strikingly simple deformation argument to show that the map
$\Om|_{H^\al_0(\T)} : H^\al_0(\T)\to\h^{\al+1/2}$ is onto. 
More precisely it uses that $(1)$, $d_0\Om_\al : H^\al_0(\T)\to\h^{\al+1/2}$ is a linear isomorphism, $(2)$, that the map
$\Om : H^{-1}_0(\T)\to\h^{-1/2}$ is a canonical bi-analytic diffeomorphism, and $(3)$, that the Hamiltonian vector field
defining the deformation is actually in $L^2$.
The same method could also be used for the proof of analogous results for more general weighted Sobolev spaces.
In a subsequent work we plan to apply our technique to the defocusing Nonlinear Schr{\"o}dinger equation.

\vspace{0.3cm}

\noindent{\em Related work:} Theorem \ref{Th:main} improves on earlier results in \cite{KT3} where it was shown that 
$\Om|_{H^\al_0(\T)} : H^\al_0(\T)\to\h^{\al+1/2}$ is a bianalytic diffeomorphism onto its image for any $-1<\al<0$.
For partial results in this direction see also \cite{Mohr}. The statement of Theorem \ref{Th:characterization} adds to
numerous results characterizing the regularity of a potential by the decay of the corresponding gap lengths --
see e.~g. \cite{DjMit}, \cite{KapMit}, \cite{Koro}, \cite{Mar}, \cite{Po} and references therein.
However only a few results concern potentials in spaces of distributions -- see \cite{KapMohr}, \cite{Koro}
(cf. also \cite{KT3} and the references therein). In a first attempt we have tried to apply the most beautiful and
most simple approach among all the papers cited, due to P{\"o}schel \cite{Po}, to our case.
However his methods seem to fail if $\al\le -3/4$. 

The idea of using flows to prove that a map is onto is not new in this subject.
It has been used e.g. by P{\"o}schel and Trubowitz in their book \cite{PT} or, to give
a more recent example, in work of Chelkak and Korotyaev \cite{ChK}. 
More precisely, in \cite[Theorem 2, p. 115]{PT}, the authors use flows to characterize
sequences coming up as as sequences of Dirichlet eigenvalues of Schr{\"o}dinger operators
$-\frac{d^2}{dx^2}+q$ on $[0,1]$ with an even $L^2$-potential $q$.
Note however, that in this paper the use of flows is of a different nature, best explained by the fact
that they are regularizing - in other words, the vector fields describing the deformations are 
in a higher Sobolev space than the underlying phase space.

\section{Proof of Theorem \ref{Th:main}}
Let $\Om$ be the Birkhoff map $\Om : H^{-1}_0(\T)\to\h^{-1/2}$ constructed in \cite{KT3}
-- see also Appendix for a brief summary of the results in \cite{KT3}.
By Theorem \ref{Th:main*} in Appendix, the Birkhoff map $\Om$ is onto and
for any given $\al>-1$ its restriction to $H^\al_0(\T)$ is a map 
\begin{equation}\label{e:Om_al}
\Om_\al:=\Om|_{H^\al_0(\T)} : H^\al_0(\T)\to\h^{\al+1/2}
\end{equation}
which is a bianalytic diffeomorphism onto its image.
Hence, in order to prove Theorem \ref{Th:main} we need to prove that \eqref{e:Om_al} is onto.

Assume that there exists $-1\le\al\le 0$ such that $\Om_\al : H^\al_0(\T)\to\h^{\al+1/2}$ is {\em not} onto.
As $\Om : H^{-1}_0(\T)\to\h^{-1/2}$ is onto it then follows that there exists 
\begin{equation}\label{e:assumption}
q_0\in H^{-1}_0(\T)\setminus H^\al_0(\T)
\end{equation}
such that $\Om(q_0)\in\h^{\al+1/2}$. 

As $\Om(0)=0$ and as by Corollary \ref{Coro:neighborhood_of_zero} in the Appendix below the differential
\[
d_0\Om_\al : H^\al_0(\T)\to\h^{\al+1/2}
\]
of \eqref{e:Om_al} at $q=0$ is a linear isomorphism, one gets from the inverse function theorem that there exist
an open neighborhood $U_\al$ of zero in $H^\al_0(\T)$ and an open neighborhood $V_\al$ of zero in $\h^{\al+1/2}$ such that 
\begin{equation}\label{e:diffeomorphism_near_zero}
\Om|_{U_\al} : U_\al\to V_\al
\end{equation}
is a diffeomorphism.

Recall that for any $k\ge 1$ the angle variable $\theta_k$ constructed in \cite{KT3} is a real-analytic function on 
$H^{-1}_0(\T)\setminus D_k$ with values in $\R/2\pi\Z$ where $D_k:=\{q\in H^{-1}_0(\T)\,|\,\ga_k(q)=0\}$ is
a real-analytic sub-variety in $H^{-1}_0(\T)$ (cf. Appendix). 
As $\theta_k$ is real-analytic, the mapping $H^{-1}_0(\T)\setminus D_k\to H^1_0(\T)$, 
$q\mapsto\frac{\partial\theta_k}{\partial q}(q)$,
is real-analytic\footnote{$\frac{\partial\theta_k}{\partial q}$ denotes the $L^2$-gradient of $\theta_k$.}
and therefore,
\begin{equation}\label{e:Y_l}
H^{-1}_0(\T)\setminus D_k\to L^2_0(\T),\,\,\,q\mapsto Y_k(q):=\frac{d}{dx}\,\frac{\partial\theta_k}{\partial q}(q)\,,
\end{equation}
is real-analytic as well. Then $Y_k$ is a Hamiltonian vector field on $H^{-1}_0(\T)\setminus D_k$, which
defines a dynamical system
\begin{equation}\label{e:ds}
{\dot q}=Y_k(q),\,\,\,q(0)=q_0\in H^{-1}_0(\T)\setminus D_k\,.
\end{equation}
Let $q_0\in H^{-1}_0(\T)\setminus D_k$ and assume that
\[
\Om(q_0)=(z_1^0,z_2^0,...)\in\h^{\al+1/2}
\]
where for any $n\ge 1$, $z_n^0=(x_n^0,y_n^0)$.
Take $\varepsilon>0$ such that the ball 
\[
B(2\varepsilon):=\{z\in\h^{\al+1/2}\,|\,\|z\|_{\al+1/2}<2\varepsilon\}
\]
is contained in the neighborhood $V_\al$ of zero in $\h^{\al+1/2}$ chosen in \eqref{e:diffeomorphism_near_zero}.
Denote by $I_n=I_n(q)$ the $n$'th action variable of a potential $q$ -- see Appendix.
Note that for any $q$ in $H^{-1}_0(\T)$
\begin{equation}\label{e:I_n}
2\,I_n(q)=\|z_n(q)\|^2=x_n(q)^2+y_n(q)^2
\end{equation}
where $\Om(q)=(z_n(q))_{n\ge 1}$ and $z_n(q)=(x_n(q),y_n(q))$.
Consider the sequence of potentials $(q_n)_{n\ge 1}$ in $H^{-1}_0(\T)$ defined recursively for $n\ge 1$ by
\[
q_n:=\left\{
\begin{array}{l}
q_{n-1}\;\;\;\mbox{if}\;\;\;2I_n(q_{n-1})<\varepsilon/(n^{1+2\al}\,2^n)\\
(q_{n-1})_{,n}\;\;\;\mbox{otherwise}
\end{array}
\right.
\]
where $(q_{n-1})_{,n}$ is obtained by shifting $q_{n-1}$ along the flow of the vector field $Y_n$ such that
\[
2I_n((q_{n-1})_{,n})<\varepsilon/(n^{1+2\al}\,2^n)\,.
\] 
The existence of $(q_{n-1})_{,n}$ follows from Lemma \ref{Lem:ds} $(a)$ below. Moreover, by the commutator relations
\eqref{e:commutators1} in Appendix,
\[
Y_n(I_m)=\{I_m,\theta_n\}=0\;\;(\;n\ne m),
\]
the vector field $Y_n$ preserves the values of the action variables $I_m$ for any $m\ne n$.
In particular, we get
\begin{equation}\label{e:ineq1}
2 I_j(q_n)\le\varepsilon/(j^{1+2\al}\,2^j),\;\;\;\forall\;1\le j\le n
\end{equation}
and
\begin{equation}\label{e:ineq2}
2 I_j(q_n)=\|z_j^0\|^2,\;\;\;\forall j>n\,.
\end{equation} 
One obtains from \eqref{e:ineq1}, \eqref{e:ineq2}, and $\|z_j\|^2=2I_j$ (cf.\ \eqref{e:z_k}) that
\begin{equation}\label{e:inequality}
\|\Om(q_n)\|_{\al+1/2}^2=\sum_{j=1}^\infty j^{1+2\al}\|z_j(q_n)\|^2\le
\varepsilon\sum\limits_{1\le j\le n}\frac{1}{2^j}+
\sum_{j\ge n+1}j^{1+2\al}\|z_j^0\|^2\,.
\end{equation}
As $\sum_{j\ge 1}j^{1+2\al}\|z_j^0\|^2=\|\Om(q_0)\|_{\al+1/2}^2<\infty$,
one gets from \eqref{e:inequality} that there exists $N\ge 1$ such that for any $n\ge N$,
$\|\Om(q_n)\|_{\al+1/2}<2\varepsilon$\,. In particular,
$\Om(q_N)\in V_\al$ and, as $\Om|_{U_\al} : U_\al\to V_\al$ is a diffeomorphism,
the bijectivity of the Birkhoff map $\Om : H^{-1}_0\to\h^{-1/2}$ implies that
\begin{equation}\label{e:q_N-regularity}
q_N\in U_\al\subseteq H^\al_0(\T)\,.
\end{equation}
On the other side, it follows from \eqref{e:assumption} and Lemma \ref{Lem:ds} $(b)$ that 
\[
(q_n)_{n\ge 1}\subseteq H^{-1}_0(\T)\setminus H^\al_0(\T)
\]
which implies $q_N\in H^{-1}_0(\T)\setminus H^\al_0(\T)$, contradicting \eqref{e:q_N-regularity}.
This completes the proof of Theorem \ref{Th:main}.
\finishproof

The following Lemma was used in the proof of Theorem \ref{Th:main}.
\begin{Lemma}\hspace{-2mm}{\bf .}\label{Lem:ds}
For any $k\ge 1$ and for any initial data $q_0\in H^{-1}_0(\T)\setminus D_k$ the initial value problem \eqref{e:ds}
has a unique solution in $C^1((-I_k^0,\infty), H^{-1}_0(\T))$ where $I_k^0\ge 0$ is the value of the action variable $I_k$
at $q_0$. The solution has the following additional properties:
\begin{itemize}
\item[$(a)$] $\lim\limits_{t\to-I^0_k+0}I_k(q(t))=0$;
\item[$(b)$] $q(t)-q_0\in L^2_0(\T)$.
\end{itemize}
\end{Lemma}
\noindent{\em Proof of Lemma \ref{Lem:ds}.}
By Theorem \ref{Th:main*} in the Appendix, the Birkhoff map $\Om : H^{-1}_0(\T)\to\h^{-1/2}$,
\[
q\mapsto\Om(q)=(z_1,z_2,...),\;\;\;z_n=(x_n,y_n),
\]
is a bianalytic diffeomorphism that transforms the Poisson structure $\frac{d}{dx}$ on $H^{-1}_0(\T)$ (cf.\ Appendix)
into the canonical Poisson structure on $\h^{-1/2}$ defined by the relations $\{x_m,x_n\}=\{y_m,y_n\}=0$ and
$\{x_m,y_n\}=\delta_{mn}$ that hold for any $m,n\ge 1$.\footnote{Here $\delta_{mn}$ denotes the Kronecker delta.}
Moreover, it follows from the construction of the Birkhoff map $\Om$ that
$\theta_k$ is the argument of the complex number $x_k+iy_k$. In particular,
in Birkhoff coordinates $(z_1,z_2,...)\in\h^{-1/2}$, one has for any $q\in H^{-1}_0(\T)\setminus D_k$
\begin{equation}\label{e:Y_l*}
d\Om(Y_k)=\frac{x_k}{x_k^2+y_k^2}\,\frac{\partial}{\partial x_k}+\frac{y_k}{x_k^2+y_k^2}\,\frac{\partial}{\partial y_k}\,.
\end{equation}
The dynamical system corresponding to the vector field \eqref{e:Y_l*} in $\h^{-1/2}$ has a unique solution for any
initial data $(x_n^0,y_n^0)_{n\ge 1}$ that is defined on the time interval $(-((x_k^0)^2+(y_k^0)^2)/2,\infty)$.
Hence, as $\Om : H^{-1}_0(\T)\to\h^{-1/2}$ is a diffeomorphism, the dynamical system \eqref{e:ds} has a unique solution $q(t)$
on $H^{-1}_0(\T)\setminus D_k$ defined for $t\in (-I_k(q_0),\infty)$.
Moreover, one gets from \eqref{e:Y_l*} and \eqref{e:I_n} that
\[
\lim\limits_{t\to-I_k(q_0)+0}I_k(q(t))=0\,.
\]
This completes the proof of $(a)$.
In order to prove $(b)$ we integrate both sides of \eqref{e:ds} in $H^{-1}_0(\T)$ and get that
for any $t\in(-\infty,I_k(q_0))$,
\begin{equation}\label{e:integral_equation}
q(t)=q_0+\int_0^tY_k(q(s))\,ds\,.
\end{equation}
As the mapping \eqref{e:Y_l} is real-analytic (and hence, continuous) and as the solution $q(t)$ of \eqref{e:ds}
is a $C^1$-curve $(-\infty,I_k(q_0))\to H^{-1}_0(\T)$, the integrand in \eqref{e:integral_equation}
is in $C^0((-I_k(q_0),\infty),L^2_0(\T))$. In particular, the integral in \eqref{e:integral_equation} converges
with respect to the $L^2$-norm, and hence represents an element in $L^2_0(\T)$.
This proves $(b)$.
\finishproof

\section{Proof of Theorem 2}
As for any constant $c \in \R,$ the potentials $q$ and $q+c$ have the same sequence of gap lengths
$(\ga_k)_{k\ge 1}$ it is enough to prove the statement of the theorem for $q\in H^{-1}_0(\T)$.

For $q\in H^{-1}_0(\T)$ given let
\[
z=(z_1,z_2,...)=\Om(q),
\]
where for any $n\ge 1$, $z_n=(x_n,y_n)$.
By Proposition \ref{Prop:actions<->gaps} in Appendix, there exist constants $0<C_1<C_2<\infty$ and $n_0\ge 1$
depending on $q$ such that for any $n\ge n_0$,
\begin{equation}\label{e:actions-gaps}
C_1\,\frac{\ga_n^2}{n}\le I_n\le C_2\,\frac{\ga_n^2}{n}
\end{equation} 
where $I_n$ is the $n$-th action variable of the
given potential $q$.
Using that
\[
I_n=(x_n^2+y_n^2)/2
\]
we get from \eqref{e:actions-gaps} that for any given $\alpha\ge -1$,
\begin{equation}\label{e:equivalence1}
(z_n)_{n\ge 1}\in\h^{\al+1/2}\;\;\;\Longleftrightarrow\;\;\;(\ga_n)_{n\ge 1}\in\h^\al\,.
\end{equation}
On the other side, it follows from Theorem \ref{Th:main} and the injectivity of
$\Om : H^{-1}_0(\T)\to\h^{-1/2}$
that
\begin{equation}\label{e:equivalence2}
(z_n)_{n\ge 1}\in\h^{\al+1/2}\;\;\;\Longleftrightarrow\;\;\;q\in H^\al_0(\T)\,.
\end{equation}
Theorem \ref{Th:characterization} now follows from
\eqref{e:equivalence1} and \eqref{e:equivalence2}.
\finishproof

\section{Appendix}
In this appendix we collect the properties of the Birkhoff map $\Om : H^{-1}_0(\T)\to\h^{-1/2}$
constructed in \cite{KT3} that were used in the proofs of Theorem \ref{Th:main} and
Theorem \ref{Th:characterization}. 

The Korteweg - de Vries equation (KdV) 
\begin{eqnarray*}
q_t-6qq_x+q_{xxx}&=&0\\
q|_{t=0}&=&q_0
\end{eqnarray*}
on the circle can be viewed as an integrable PDE, i.e.\ an integrable Hamiltonian system of   
infinite dimension. As a phase space we consider the Sobolev space $H^\al(\T)$ $(\al\ge -1)$ 
of real valued distributions on the circle. The Poisson bracket is the one proposed by Gardner,   
\begin{equation}\label{e:PB}
\{F,G\}:=\int_{\T}\frac{\partial F}{\partial q}\frac{d}{dx}\Big(\frac{\partial G}{\partial q}\Big)\;dx
\end{equation}
where $F$, $G$ are $C^1$-functions on $H^\al(\T)$ and $\frac{\partial F}{\partial q}$,
$\frac{\partial G}{\partial q}$ denote the $L^2$-gradients of $F$ and $G$ respectively which
are assumed to be sufficiently smooth so that the Poisson bracket is well defined.
For $q$ sufficiently smooth, i.e. $q\in H^1_0(\T)$, the Hamiltonian ${\mathcal H}$ corresponding to KdV is given by  
\[
{\mathcal H}(q)=\int_{\T}((\partial_xq)^2/2+q^3)\;dx 
\]
and the KdV equation can be written in Hamiltonian form
\[
q_t=\frac{d}{dx}\,\frac{\partial{\mathcal H}}{\partial q}\,.
\]
Note that the Poisson structure is degenerate and admits the average $[q]:=\int_{\T}q(x)\,dx$ as
a Casimir function. Moreover, the Poisson structure is regular and induces a trivial foliation whose leaves   
are given by
\[   
H^{\al}_c(\T)=\{q\in H^\al(\T)\;|\;[q]=c\}\,.   
\]   
Introduce the set
\[
D_k:=\{q\in H^{-1}_0(\T)\,|\,\ga_k(q)=0\}\,.
\]
For any $q\in H^{-1}_0(\T)\setminus D_k$ define
\begin{equation}\label{e:z_k}
z_k(q):=\sqrt{2I_k(q)}\,\Big(\cos(\theta_k(q)),\sin(\theta_k(q))\Big)\,,
\end{equation}
where $I_k(q)$ is the $k$'th action variable and $\theta_k(q)$ is the $k$'th angle variable
of the KdV equation (cf.\ \S\;3,\,4\;in\;\cite{KT3}). 
It is shown in \cite[\S\,5]{KT3} that the mapping $H^{-1}_0\setminus D_k\to{\mathbb R}^2$, $q\mapsto z_k(q)$, extends analytically to $H^{-1}_0(\T)$.
For any $q\in H^{-1}_0(\T)$ the action variables $(I_k)_{k\ge 1}$ of KdV are defined in terms of the periodic spectrum of
the Schr{\"o}dinger operator $-\frac{d^2}{dx^2}+q$ using the same formulas as in \cite{FM} (cf. also \cite{KP}).
For any given $\al\ge-1$ and for any $k\ge 1$ the action $I_k$ is a real analytic function on $H^{\al}_0(\T)$ (cf.\ Proposition 3.3 in \cite{KT3}).
The angle $\theta_k$ is defined modulo $2\pi$ and is a real analytic function on $H^{\al}_0(\T)\setminus (D_k\cap H^\al_0)$, where
$D_k\cap H^\al_0=\{q\in H^{\al}_0(\T)\,|\,\ga_k(q)=0\}$ is a real analytic sub-variety in $H^{\al}_0(\T)$ of
co-dimension two (cf.\ Proposition 4.3 in \cite{KT3}). By \S\;6 in \cite{KT3} we have the following
commutator relations
\begin{equation}\label{e:commutators0}
\{I_m,I_n\}=0\;\;\mbox{on}\;\; H^{-1}_0(\T)
\end{equation}
\begin{equation}\label{e:commutators1}
\{I_m,\theta_n\}=\delta_{nm}\;\;\mbox{on}\;\; H^{-1}_0(\T)\setminus D_n
\end{equation}
and
\begin{equation}\label{e:commutators2}
\{\theta_m,\theta_n\}=0\;\;\mbox{on}\;\;H^{-1}_0(\T)\setminus(D_m\cup D_n)
\end{equation}
for any $m,n\ge 1$. 
For any $q\in H^{-1}_0(\T)$ define
\begin{equation*}
\Om(q):=(z_1(q),z_2(q),...)
\end{equation*}
where $z_k=z_k(q)$ is given by \eqref{e:z_k}. It is shown in \cite{KT3} that $\Om(q)\in\h^{-1/2}$.
Recall that, for any $\al\in\R$, $\h^\al$ denotes the Hilbert space
\[
\h^\al=\{z=(x_k,y_k)_{k\ge 1}\;|\;\|z\|_\al<\infty\},
\]
with the norm
\[
\|z\|_\al:=\Big(\sum_{k\ge 1}k^{2\al}(x_k^2+y_k^2)\Big)^{1/2}\,.
\]
We supply $\h^{-1/2}$ with a Poisson structure defined by the relations
$\{x_m,x_n\}=\{y_m,y_n\}=0$ and $\{x_m,y_n\}=\delta_{mn}$ valid for any $m,n\ge 1$.
The following result is proved in \cite{KT3}.
\begin{Th}\hspace{-2mm}{\bf .}\label{Th:main*}
The mapping $\Om : H^{-1}_0(\T)\to\h^{-1/2}$ satisfies the following properties:
\begin{itemize}
\item[(i)] $\Om$ is a bianalytic diffeomorphism that preserves the Poisson bracket;
\item[(ii)] for any $\al>-1$, the restriction $\Om_\al\equiv\Om|_{H^\al_0(\T)}$ is a map 
$\Om|_{H^\al_0(\T)} : H^\al_0(\T)\to\h^{\al+1/2}$ which is one-to-one and bianalytic onto its image.
In particular, the image is an open subset in $\h^{\al+1/2}$.
\end{itemize}
\end{Th}
\begin{Coro}\hspace{-2mm}{\bf .}\label{Coro:neighborhood_of_zero}
For any $\al>-1$,
\[
d_0\Om_\al : H^{\al}_0(\T)\to\h^{\al+1/2}\,,
\]
is a linear isomorphism.
\end{Coro}
We will also need the following Proposition (cf.\ \cite[\S\,3]{KT3}).
\begin{Prop}\hspace{-2mm}{\bf .}\label{Prop:actions<->gaps}
There exists a complex neighborhood ${\mathcal W}$ of $H^{-1}_0(\T)$ in the complex space $H^{-1}_0(\T,\C)$
such that the quotient $I_n/\ga_n^2$, defined on $H^{-1}_0(\T)\setminus D_n$, extends analytically to ${\mathcal W}$ for all $n$.
Moreover, for any $\varepsilon>0$ and any $p\in{\mathcal W}$ there exists $n_0\ge 1$ and an open neighborhood
$U(p)$ of $p$ in ${\mathcal W}$ so that 
\[
\left|8\pi n\,\frac{I_n}{\ga_n^2}-1\right|\le\varepsilon
\]
for any $n\ge n_0$ and for any $q\in U(p)$.
\end{Prop}
Further we recall that for any $q\in H^{-1}_0(\T)$ one has that {\em $I_n(q)=0$ if and only if $\ga_n(q)=0$}. 
In particular, one concludes from \eqref{e:z_k} and the fact $\ga_n(0)=0 \, \forall n \ge 1$ that $\Om(0)=0$.

\end{document}